\documentclass[a4paper,11pt]{amsart}
\usepackage{amssymb,amsfonts,amsxtra,amscd}
\usepackage[all]{xy}
\usepackage{fullpage}
\theoremstyle{plain}
\newtheorem{theorem}{Theorem}[section]
\newtheorem{lemma}[theorem]{Lemma}

\newtheorem{prop}[theorem]{Proposition}
\theoremstyle{definition}
\newtheorem{defi}[theorem]{Definition}

\theoremstyle{remark}
\newtheorem{rem}[theorem]{Remark}
\numberwithin{equation}{section}
\newcommand{\ci}{\ensuremath{C_\infty}}
\newcommand{\ai}{\ensuremath{A_\infty}}

\newcommand{\ctalg}[1]{\ensuremath{\widehat{T}\Sigma #1^*}}
\newcommand{\clalg}[1]{\ensuremath{\widehat{L}\Sigma #1^*}}

\newcommand{\deof}[2][]{\ensuremath{\Omega^1_\mathrm{#1}(#2)}}
\newcommand{\drof}[2][]{\ensuremath{DR^1_\mathrm{#1}(#2)}}

\newcommand{\drzf}[2][]{\ensuremath{DR^0_\mathrm{#1}(#2)}}

\newcommand{\de}[2][]{\ensuremath{\Omega^\bullet_\mathrm{#1}(#2)}}

\newcommand{\Cbr}[2][\bullet]{\ensuremath{C_{#1}^\mathrm{Bar}(#2)}}
\newcommand{\Choch}[3][\bullet]{\ensuremath{C_{#1}^\mathrm{Hoch}(#2,#3)}}
\newcommand{\Hhoch}[3][\bullet]{\ensuremath{H_{#1}^\mathrm{Hoch}(#2,#3)}}

\newcommand{\choch}[3][\bullet]{\ensuremath{C^{#1}_\mathrm{Hoch}(#2,#3)}}
\newcommand{\hhoch}[3][\bullet]{\ensuremath{H^{#1}_\mathrm{Hoch}(#2,#3)}}

\newcommand{\cbr}[2][\bullet]{\ensuremath{C^{#1}_\mathrm{Bar}(#2)}}

\newcommand{\cchoch}[2][\bullet]{\ensuremath{CC^{#1}_\mathrm{Hoch}(#2)}}
\newcommand{\hchoch}[2][\bullet]{\ensuremath{HC^{#1}_\mathrm{Hoch}(#2)}}

\newcommand{\cotimes}{\ensuremath{\hat{\otimes}}}
\newcommand{\gf}{\ensuremath{\mathbb{Q}}}

\newcommand{\noproof}{\begin{flushright} \ensuremath{\square} \end{flushright}}
\DeclareMathOperator{\Map}{Map}
\DeclareMathOperator{\Hom}{Hom}
\DeclareMathOperator{\Der}{Der}

\DeclareMathOperator{\ad}{ad}

\begin{document}
\begin{abstract}
In this paper we  establish the existence of certain structures on the ordinary and equivariant homology of the free loop space on a manifold or, more generally, a formal Poincar\'e duality space. These structures; namely the loop product, the loop bracket and the string bracket, were introduced and studied by Chas and Sullivan under the general heading `string topology'. Our method is based on obstruction theory for $C_\infty$-algebras and rational homotopy theory. The resulting string topology operations are manifestly homotopy invariant.
\end{abstract}
\title{Symplectic \ai-algebras and string topology operations}
\author{Alastair Hamilton \and Andrey Lazarev}

\maketitle

\section{Introduction}
In this section we apply the  results  of \cite{HL1}, \cite{HL2} to establish the existence of certain structures on the ordinary and equivariant homology of the free loop space on a manifold or, more generally, a formal Poincar\'e duality space. These structures; namely the loop product, the loop bracket and the string bracket, were introduced and studied in \cite{chasul} under the general heading `string topology'. Since the original work of Sullivan and Chas a number of other papers appeared which investigated string topology structures from different perspectives. In particular, the approach of \cite{cohenj} and  \cite{jklein} was purely homotopy theoretic while that of \cite{fethvi} was based on the theory of minimal models. Ideas similar to those used in the present paper were emplyed in \cite{laz} to establish the existence of a `string bracket' on the $S^1$-equivariant homology of a punctured closed manifold (or a space homotopy equivalent to it).

It needs to be stressed that the we define `string topology operations' in a way different from that of \cite{chasul}. In view of the similarities of their properties to the operations defined by Sullivan and Chas it seems very likely that our definition is, in fact, equivalent to theirs. However, this issue is not discussed in the present paper. In the recent work \cite{CKS} the authors the authors prove that the original Sullivan-Chas operations are also homotopy invariant.

This paper is organized as follows. In Section \ref{recall} we recall the definition of a symplectic $A_\infty$ and $C_\infty$ algebra and the main result of \cite{HL2} which states, roughly speaking, that given a $C_\infty$-algebra with a compatible structure of a Frobenius algebra on its cohomology, there is a unique up to homotopy way to extend it to a symplectic $C_\infty$-algebra.

In Section \ref{strings} we use the afore-mentioned result to define the analogues of the Sullivan-Chas string topology operations on the free loop space of a simply-connected Poincar\'e duality space and on the $S^1$-equivariant version thereof. The homotopy invariance of our constructions is immediate.
\subsection{Notation and conventions}
Throughout the paper we we will assume that our ground ring is the field of rational numbers $\mathbb{Q}$; all the unmarked tensors are understood to be taken over $\gf$.

Given a graded vector space $V$ we define the tensor algebra $TV$ by
\[ TV:=\gf \oplus V \oplus V^{\otimes 2} \oplus \ldots \oplus V^{\otimes n} \oplus \ldots. \]
We define the free Lie algebra $LV$ as the Lie subalgebra of the commutator algebra $TV$ which consists of linear combinations of Lie monomials.

We will use $~\widehat{ \ }~$ to denote completion. Given a finite-dimensional vector space $V$ we can define the completed versions of $TV$ and $LV$. These will be called the \emph{pro-free} associative and Lie algebra respectively. For instance the pro-free associative algebra would be
\[ \widehat{T}V :=\prod_{i=0}^{\infty} V^{\otimes i} \]
 The pro-free Lie algebra is the Lie subalgebra of $\widehat{T}V$ consisting of all convergent infinite linear combinations of Lie monomials.

Given a pro-free associative or Lie algebra $X$, the Lie algebra consisting of all \emph{continuous} derivations
\[ \xi : X \to X \]
is denoted by $\Der(X)$. Continuous derivations on $X$ will also be called vector fields on $X$ to emphasise similarity with the corresponding geometric notion.

Our convention will be to always work with cohomologically graded objects and we consequently define the suspension of a graded module $V$ as $\Sigma V$ where $\Sigma V^i:=V^{i+1}$. We define the desuspension of $V$ as $\Sigma^{-1}V$ where $\Sigma^{-1}V^{i}:= V^{i-1}$. The term `differential graded algebra' will be abbreviated as `DGA'.

Given a graded module $V$, we denote the graded \gf-linear dual
\[ \Hom_\gf(V,\gf) \]
by $V^*$.
For the sake of clarity, when we write $\Sigma V^*$ we mean the graded module $\Hom_\gf(\Sigma V,\gf)$.

For a topological space $X$ (always assumed to be a $CW$-complex) we denote by $C^\bullet(X)$ the singular cochain algebra of $X$ with coefficients in $\mathbb{Q}$. The homology and cohomology of $X$ with coefficients in $\mathbb{Q}$ will be denoted by $H_\bullet(X)$ and $H^\bullet(X)$ respectively.   The space of all maps $S^1\rightarrow X$ from the circle $S^1$ to $X$ will be denoted by $LX$. The space $LX$ has an action of $S^1$ and the corresponding equivariant homology will be denoted by $H^{S^1}_\bullet(LM)$:
\[  H^{S^1}_\bullet(LM):= H^\bullet(ES^1\times _{S^1} LM).\]

\section{Symplectic $C_\infty$-algebras and their obstruction theory}\label{recall}
In this section the notion of $A_\infty-$ and $C_\infty$-algebra will be introduced as well as their symplectic versions and theoir obstruction theory.  This material is borrowed from \cite{HL1, HL2} and the reader is advised to refer to these papers for further details and proofs.

The concept of an $A_\infty$-algebra was devised by Stasheff \cite{staha1, staha2}; it is a homotopy invariant version of a strictly associative algebra whilst a $C_\infty$-algebra  is designed to represent homotopy commutative structures.

The notion of a \emph{symplectic} $A_\infty$-algebra (also known as an $A_\infty$-algebra with an invariant inner product or a cyclic $A_\infty$-algebra) is due to Kontsevich, \cite{kontsg, kontfd}; it could be readily extended to the context of $C_\infty$-algebras. Symplectic $A_\infty$ (respectively $C_\infty$-algebras) are homotopy invariant versions of Frobenius algebras (noncommutative and commutative respectively).

Let $V$ be a vector space which we will assume to be finite-dimensional for simplicity. Choose a basis $t_1,\ldots t_n$ of $\Sigma V^*$, then $\ctalg{V}=\gf\langle\langle t_1,\ldots, t_n\rangle\rangle$. Any continuous derivation $\xi:\ctalg{V}\to\ctalg{V}$ could be written as
\[ \xi=\sum_{i=1}^nA_i(t_1,\ldots, t_n)\partial_{t_i}. \] where $A_i(t_1,\ldots, t_n)$ are noncommutative formal power series in $t_1,\ldots, t_n$. If all the $A_i$'s have degree $i$ we say that $\xi$ has degree $i$. Thus, any such $\xi$ could be written as $\xi=\xi_0+\xi_1+\ldots$ where $\xi_0 $ is a constant derivation, $\xi_1$ is linear, $\xi_2$ is quadratic etc.

We say that a continuous derivation.  $\xi:\ctalg{V}\to\ctalg{V}$ vanishes at zero if it has the form
\[ \xi=\sum_{i \in I}A_i(t_1,\ldots, t_n)\partial_{t_i}, \]
where the power series $A_i(t_1,\ldots, t_n)$ have vanishing constant terms.

We will now recall the definition of an  {\ai} and \ci-structure on a graded vector space $V$. An  {\ai} or \ci-algebra is a graded vector space together with an  {\ai} or \ci-structure. 
\begin{defi} \label{def_infstr}
Let $V$ be a graded vector space:
\begin{enumerate}
\item[(a)]
An \ai-structure on $V$ is a vector field
\[ m:\ctalg{V} \to \ctalg{V} \]
of degree one and vanishing at zero, such that $m^2=0$.
\item[(b)]
A \ci-structure on $V$ is a vector field
\[ m:\clalg{V} \to \clalg{V} \]
of degree one, such that $m^2=0$.
\end{enumerate}
An {\ai} or \ci-algebra structure is called \emph{minimal} if its linear term $m_1$ vanishes. Note that $m_1$ determines a differential of the graded vector space  $\Sigma V^*$.
\end{defi}

\begin{defi} \label{def_infmor}
Let $V$ and $U$ be graded vector spaces:
\begin{enumerate}
\item[(a)]
Let $m$ and $m'$ be \ai-structures on $V$ and $U$ respectively. An \ai-morphism from $V$ to $U$ is a continuous algebra homomorphism
\[ \phi:\ctalg{U} \to \ctalg{V} \]
of degree zero such that $\phi \circ m'=m \circ \phi$.
\item[(b)]
Let $m$ and $m'$ be \ci-structures on $V$ and $U$ respectively. A \ci-morphism from $V$ to $U$ is a continuous algebra homomorphism
\[ \phi:\clalg{U} \to \clalg{V} \]
of degree zero such that $\phi \circ m'=m \circ \phi$.
\end{enumerate}
If the linear part $\phi_1$ of $\phi$ determines a quasi-isomorphism $\Sigma U^*\rightarrow \Sigma V^*$ then $\phi$ is called an $A_\infty-$ (respectively $C_\infty-$) isomorphism.
\end{defi}

\begin{rem} \label{rem_calseq}
Any continuous derivation $m:\clalg{V} \to \clalg{V}$ can be uniquely extended to a continuous Hopf algebra derivation (where \ctalg{V} is equipped with the shuffle coproduct)
\[ m : \ctalg{V} \to \ctalg{V}; \]
clearly all continuous Hopf algebra derivations on \ctalg{V} are obtained from $\Der(\clalg{V})$ in this manner. It follows that any \ci-structure gives rise to an \ai-structure in this way
and all continuous Hopf algebra homomorphisms are obtained from continuous Lie algebra homomorphisms in this way. This is because the Lie subalgebra of primitive elements of our Hopf algebra \ctalg{V} coincides with the Lie subalgebra \clalg{V}.

The category $\ai\mathrm{-}algebras$ has a subcategory whose objects consist of the \ai-algebras whose \ai-structure is also a Hopf algebra derivation and whose morphisms are the \ai-morphisms which are also Hopf algebra homomorphisms. It follows from the discussion above that this category is isomorphic to $\ci\mathrm{-}algebras$.
\end{rem}

\subsection{Cyclic and Hochschild cohomology of $A_\infty$-algebras}
We start by recalling the notion of noncommutative forms on a pro-free algebra.
Let $A=\widehat{T}V$ be a pro-free algebra on a vector space $V$.
\begin{defi}\
\begin{enumerate}\item
The module of noncommutative differentials \deof[]{A} is defined as
\[ \deof[]{A}:= A \cotimes (A/\gf). \]
Let us write $x\cotimes y$ as $x\cotimes dy$. \deof[]{A} has the structure of an $A$-bimodule via the actions
\begin{align*}
& a \cdot x\cotimes dy:= ax \cotimes dy, \\
& x \cotimes dy \cdot a:= x \cotimes d(ya) - xy \cotimes da.
\end{align*}
\item
 The module of noncommutative forms \de[]{A} is defined as
\[ \de[]{A}:=\widehat{T}_A(\Sigma^{-1}\deof[]{A})=A \times \prod_{i=1}^\infty \underbrace{\Sigma^{-1}\deof[]{A} \underset{A}{\cotimes} \ldots \underset{A}{\cotimes} \Sigma^{-1}\deof[]{A}}_{i \text{ factors}}. \]
Since \deof[]{A} is an $A$-bimodule, \de[]{A} has the structure of a formal associative algebra whose multiplication is the standard associative multiplication on the tensor algebra $\widehat{T}_A(\Sigma^{-1}\deof[]{A})$. The map $d:A \to \deof[]{A}$ lifts uniquely to give \de[]{A} the structure of a formal DGA.\end{enumerate}
\end{defi}
Further, the noncommutative de Rham complex  $DR^\bullet(A)$ is defined as the quotient of $\Omega^\bullet{A}$ by the subspace spanned over $\gf$ by graded commutators i.e. all elements of the form $\omega_1\cdot\omega_2-(-1)^{|\omega_1||\omega_2|}\omega_2\cdot \omega_1$ where  $\omega_i\in \Omega^\bullet, i=1,2$.

Now let $\xi$ be a vector field on $A$. Associated to $\xi$ is the operator $L_\xi$ on \deof[]{A} defined by the following axioms:
\begin{displaymath}
\begin{array}{lc}
L_\xi(x):=\xi(x), & x \in X; \\
L_\xi(dx):=(-1)^{|\xi|}d(\xi(x)), & x \in X. \\
\end{array}
\end{displaymath}
It is clear that $\xi$ determines a well-defined operator on $DR^\bullet(A)$ which will again be denoted by the symbol $L_m$.
\begin{defi}\
\begin{enumerate}\item
Let $V$ be a graded vector space:

Let $m:\ctalg{V} \to \ctalg{V}$ be an \ai-structure. The Hochschild complex of the \ai-algebra $V$ with coefficients in $V$ is defined on the module consisting of all continuous derivations of \ctalg{V}:
\[ \choch{V}{V}:=\Sigma^{-1}\Der(\ctalg{V}). \]
The differential $d:\choch{V}{V} \to \choch{V}{V}$ is given by
\[ d(\xi):=[m,\xi], \quad \xi \in \Der(\ctalg{V}). \]
The Hochschild cohomology of $V$ with coefficients in $V$ is defined as the cohomology of the complex \choch{V}{V} and denoted by \hhoch{V}{V}.
\item
Let $m:\ctalg{V} \to \ctalg{V}$ be an \ai-structure. The Hochschild complex of the \ai-algebra $V$ with coefficients in $V^*$ is defined on the module consisting of all 1-forms:
\[ \choch{V}{V^*}:=\Sigma\drof[]{\ctalg{V}}. \]
The differential on this complex is the (suspension of the) Lie operator of the vector field $m$;
\[ L_m:\drof[]{\ctalg{V}} \to \drof[]{\ctalg{V}}. \]
The Hochschild cohomology of $V$ with coefficients in $V^*$ is defined as the cohomology of the complex \choch{V}{V^*} and denoted by \hhoch{V}{V^*}.\end{enumerate}
\end{defi}
\begin{rem}
Note that  $\Der(\ctalg{V})\cong \Hom(V,\ctalg{V})$ and $DR^1(A)\cong V\cotimes \widehat{T}V$. From this it is not hard to deduce that if $(V,m)$ is in fact a strictly associative algebra (i.e. all higher products $m_i, i>2$ vanish) then the definition of \choch{V}{V^*} and of \choch{V}{V} agrees with the standard one.
\end{rem}
Here is a definition of the cyclic cohomology of an $A_\infty$-algebra.
\begin{defi}
Let $m:\ctalg{V} \to \ctalg{V}$ be an \ai-structure. The cyclic Hochschild complex of the \ai-algebra $V$ is defined on the module of 0-forms vanishing at zero:
\[ \cchoch{V}:=\Sigma\left[\{q \in \drzf[]{\ctalg{V}}:\ q \text{ vanishes at zero} \}\right]. \]
The differential on this complex is the restriction of the (suspension of the) Lie operator of the vector field $m$;
\[ L_m:\drzf[]{\ctalg{V}} \to \drzf[]{\ctalg{V}} \]
 The cyclic Hochschild cohomology of $V$ is defined as the cohomology of the complex \cchoch{V} and is denoted by \hchoch{V}.
\end{defi}

\subsection{Symplectic structures}
Let $(\widehat{T}\Sigma V^*,m)$ be an $A_\infty$-algebra. We assume that $V$ is finite dimensional over $\gf$ and that
$V$ possesses a nondegenerate graded symmetric scalar product $\langle,\rangle$ which we will refer to as the \emph{inner product}. Then $\Sigma V$ also acquires a scalar product which we will denote by the same symbol $\langle,\rangle$; namely:
\[\langle\Sigma a,\Sigma b\rangle:=(-1)^{|a|}\langle a,b\rangle.\]
It is easy to check that the product on $\Sigma V$ will be graded skew-symmetric, in other words it will determine a (linear) graded symplectic structure on $\Sigma V$.

We will consider the inner product on the underlying space of a symplectic \ai-algebra $(V,m)$ $\langle,\rangle$ as an element $\omega\in(\Sigma V^*)^{\otimes 2}$. Given a basis $x_i$ in $\Sigma V^*$ the element $\omega\in T^2(\Sigma V^*)$ could be written as $\omega=\sum\omega^{ij}x_i\otimes x_j$. Consider the element \[[\omega]:=\sum\omega^{ij}[x_i,x_j]\in L^2(\Sigma V^*)\hookrightarrow T^2(\Sigma V^*).\] Clearly $[\omega]$ does not depend on the choice of a basis. A derivation $\xi\in\Der(\ctalg{V})$ will be called \emph{symplectic} if $\xi[w]=0$.

\begin{defi}
An $A_\infty$-algebra $(\widehat{T}\Sigma V^*, m)$ with an inner product $\langle,\rangle$ is called \emph{symplectic} if ${m}$ is a symplectic derivation.
A $C_\infty$-algebra is called symplectic if it is so considered as an $A_\infty$-algebra.
\end{defi}
Let $(\widehat{T}\Sigma V^*, m_V)$ and $(\widehat{T}\Sigma W^*, m_W)$ be two symplectic $A_\infty$-algebras and $[\omega]_V, [\omega]_W$ be the corresponding elements defined above.  We say that an $A_\infty$-morphism between two symplectic $A_\infty$-algebras $f:(\widehat{T}\Sigma V^*, m_V)\rightarrow (\widehat{T}\Sigma W^*, m_W)$ is \emph{symplectic} if $f$ takes $[\omega]_V$ to $[\omega]_W$. A symplectic $C_\infty$-morphism is defined analogously.

We see, therefore, that a symplectic $C_\infty$-algebra is a homotopy associative and commutative algebra with an inner product compatible with all structure maps. In particular, the homology of a symplectic $C_\infty$-algebra is a usual commutative Frobenius algebra. The following result, proved in \cite{HL2} shows that, rather surprisingly, any minimal $C_\infty$-algebra which is symplectic \emph{only up to homotopy} could be uniquely extended to a genuine symplectic $C_\infty$-algebra.
\begin{theorem} \label{thm_main}
Let $A:=(V,m=(m_2, m_3,\ldots)$ be a minimal $C_\infty$-algebra and $\langle, \rangle$ be an inner product on $V$ which is compatible with $m_2$ in the sense that for any $a,b,c\in V$ one has $\langle ab,c\rangle=\langle a, bc\rangle$. Then there exists a minimal symplectic $C_\infty$-structure on $V$: $(V,m^\prime=m_2^\prime,m_3^\prime,\ldots)$ such that
\begin{enumerate}
\item $m_2=m_2^\prime$ and
the $C_\infty$-algebras $(V,m)$ and $(V,m^\prime)$ are $C_\infty$-isomorphic. We will call $(V,m^\prime)$ satisfying these conditions a \emph{symplectic lifting} of the $(V,m)$.
\item Any two symplectic liftings are symplectically $C_\infty$-isomorphic.
\end{enumerate}
\end{theorem}
%\subsection{Hochschild and cyclic cohomology of symplectic $A_\infty$-algebras}
If an $A_\infty$-algebra $(V,m)$ has a symplectic structure then its Hochschild cohomology with coefficients in $V$ and in $V^*$ are isomorphic:
\begin{prop} \label{lem_duaiso}
Let $(V,m,\omega)$ be

\end{prop}\noproof
We will need the following result.
\begin{prop} \label{lem_duciso}
Let $(V,m,\omega)$ be
a \emph{symplectic} \ai-algebra
then the Lie subalgebra of
$\choch{V}{V}:=\left(\Der(\ctalg{V}),\ad m\right)$ %or
consisting of all \emph{symplectic} vector fields forms a subcomplex denoted by
$S\choch{V}{V}$
respectively. Furthermore, there is an isomorphism
$\Sigma^{|\omega|-2}\cchoch{V} \to S\choch{V}{V}$ %or
\end{prop}\noproof

\section{Main construction} \label{strings}
\begin{theorem} \label{loopproduct}
Let $M$ be a simply-connected  Poincar\'e duality space of formal dimension $d$. Then the graded vector space $\mathbb{H}_\bullet(M):=H_{\bullet+d}(LM)$ has the structure of a graded Gerstenhaber algebra. Two homotopy equivalent Poincar\'e duality spaces give rise to isomorphic Gerstenhaber algebras.
\end{theorem}

\begin{rem}
The graded Gerstenhaber algebra structure on $H_{\bullet+d}(LM)$ consists of two operations: a graded commutative and associative product
\[a\otimes b\mapsto a\bullet b:H_p(LM)\otimes H_q(LM)\rightarrow H_{p+q-d}(LM)\]
and the graded Lie bracket
\[a\otimes b\mapsto \{a, b\}:H_p(LM)\otimes H_q(LM)\rightarrow H_{p+q-d+1}(LM).\]
Moreover, $\{a,?\}$ is a graded derivation of $\bullet$ for any $a\in H_\bullet(LM)$. The operations $\bullet$ and $\{,\}$ are called the loop product and the loop bracket respectively.
\end{rem}

\begin{proof}
There exists a minimal $C_\infty$-algebra $(H^\bullet(M),m)$ weakly equivalent to the cochain algebra $C^\bullet(M)$. Moreover, by Theorem \ref{thm_main} we could assume that the Poincar\'e duality form extends to an invariant inner product of order $d$ on $(H^\bullet(M),m)$. To ease notation we will denote $H^\bullet(M)$ by $V$. Then \hhoch{V}{V} will denote the Hochschild cohomology of $(H^\bullet(M),m)$ with coefficients in itself and \hhoch{V}{V^*} the Hochschild cohomology of $(H^\bullet(M),m)$ with coefficients in $V^*$.

Then by Proposition \ref{lem_duaiso} we have an isomorphism $\hhoch{V}{V^*} \cong \hhoch[\bullet+d]{V}{V}$. It is well known that the Hochschild cohomology of a differential graded algebra with coefficients in itself has the structure of a graded Gerstenhaber algebra; moreover two weakly equivalent DGAs give rise to isomorphic Gerstenhaber algebras, cf. \cite{femeth}. Since there exists a DGA weakly equivalent to the \ai-algebra $V$ we conclude that $\hhoch{V}{V} \cong \hhoch[\bullet-d]{V}{V^*}$ supports the structure of a graded Gerstenhaber algebra.

Since $C^\bullet(M)$ is weakly equivalent to $(V,m)=(H^\bullet(M),m)$ we have an isomorphism of graded vector spaces $\hhoch{C^\bullet(M)}{[C^\bullet(M)]^*} \cong \hhoch{V}{V^*}$. Therefore  \hhoch[\bullet-d]{C^\bullet(M)}{[C^\bullet(M)]^*} has the structure of a graded Gerstenhaber algebra.

Furthermore, the graded space \hhoch{C^\bullet(M)}{[C^\bullet(M)]^*} is isomorphic to the $\gf$-vector dual to \Hhoch{C^\bullet(M)}{C^\bullet(M)} (Hochschild homology of $C^\bullet(M)$) and the latter is isomorphic to $H^\bullet(LM)$ by \cite{burfi}, \cite{goodwi}.
Since the space $LM$ is of finite type we have  $[H^\bullet(LM)]^*\cong H_\bullet(LM)$ and we conclude that  \[\mathbb{H}_\bullet(M)=H_{\bullet+d}(LM)\cong \hhoch[\bullet-d]{C^\bullet(M)}{[C^\bullet(M)]^*}\]
has the structure of a graded Gerstenhaber algebra as claimed.

Homotopy invariance is likewise clear since for two homotopy equivalent spaces $M$ and $N$ the cochain algebras $C^\bullet(M)$ and $C^\bullet(N)$ are weakly equivalent.
\end{proof}

\begin{rem}
Although two homotopy equivalent Poincar\'e duality spaces $M, M^\prime$ give rise to isomorphic Gerstenhaber algebras on their string homology, this structure is not natural in the sense that a map $M\rightarrow M^\prime$ does not lead to a map between the corresponding Gerstenhaber algebras. An analogous situation arises when one considers the monoid of self-maps $\Map(X,X)$ of a topological space $X$; the association $X\mapsto \Map(X,X)$ is not a functor even though homotopy equivalent spaces give rise to homotopy equivalent monoids of self-maps. The same remark applies to the string bracket considered below.
\end{rem}

The other part of the string topology operations is called the \emph{string bracket} and is defined on the equivariant homology of $LM$: $H_\bullet^{S^1}(LM):=LM\times_{S^1}ES^{1}$. To put the string bracket in the proper context we need to introduce \emph{negative} cyclic cohomology of $A_\infty$-algebras.

Let $(V,M)$ be an $A_\infty$-algebra, not necessarily unital.

\begin{defi}
The negative cyclic complex $CC^\bullet_-(V)$ of $(V,m)$ is the total complex of the following bicomplex formed by taking direct sums:
\begin{equation}\label{negativecyclic}
\xymatrix{\ldots \ar^-{1-z}[r] & \cbr{V} \ar^-N[r] & \choch{V}{V^*} \ar^-{1-z}[r] & \cbr{V} }.
\end{equation}
The cohomology of $CC^\bullet_-(V)$ will be denoted by $HC^\bullet_-(V)$.
\end{defi}

\begin{rem}
Standard spectral sequence arguments show that two weakly equivalent \ai-algebras have isomorphic negative cyclic cohomology.
\end{rem}

\begin{lemma}\label{comparison}
Let $(V,m)$ be a unital $A_\infty$-algebra for which there exists an integer $N$ such that $\hhoch[k]{V}{V^*}=0$ for $k>N$. Then for any integer $n$ we have
\[HC^n_-(V)\cong \hchoch[n+1]{V}.\]
\end{lemma}

\begin{proof}
Note that the even-numbered columns isomorphic to \cbr{V} are acyclic since the $A_\infty$-algebra $(V,m)$ is unital. This together with
the assumption that $\hhoch[k]{V}{V^*}=0$ for sufficiently large $k$ ensures that $CC^\bullet_-(V)$ is quasi-isomorphic to the subcomplex $\overline{CC}^\bullet_-(V)$ whose columns are appropriate truncations of the columns of $CC^\bullet_-(V)$  and which has no nonzero terms above a certain horizontal line.

Consider the auxiliary complex $\widetilde{CC}^\bullet_-(V)$ formed by the \emph{direct product} totalisation of \eqref{negativecyclic}
Then the complex  $\overline{CC}^\bullet_-(V)$ can also be considered as a subcomplex of $\widetilde{CC}^\bullet_-(V)$. Note that for the complex $\overline{CC}^\bullet_-(V)$ there is no difference between the direct product or direct sum totalisation and therefore both its spectral sequences converge.

Comparing the appropriate spectral sequences for $\overline{CC}^\bullet_-(V)$ and $\widetilde{CC}^\bullet_-(V)$ we see, using the exactness of  the rows, that these complexes are quasi-isomorphic and therefore they are both quasi-isomorphic to $CC^\bullet_-(V)$.
That shows that under our assumptions the direct product and direct sum totalisations of \eqref{negativecyclic} have the same cohomology.

Next consider the following bicomplex $\widehat{CC}^\bullet_-(V)$:
\[\xymatrix{ \ldots \ar^-{N}[r] & \choch{V}{V^*} \ar^-{1-z}[r] & \cbr{V} \ar^-{N}[r] & \choch{V}{V^*} }.\]
Its cohomology will be denoted by $\widehat{HC}^\bullet_-(V)$. Similarly to the above, it makes no difference which totalisation (direct sum or direct product) we take.
Now zig-zag arguments (or the appropriate spectral sequence)  show, using the exactness of the rows,  that  $\widehat{HC}^\bullet_-(V)$ is isomorphic to the complex formed by the image of the horizontal differential in the rightmost column. The latter complex is clearly isomorphic to \cchoch{V}. Therefore   $\widehat{HC}^\bullet_-(V)\cong HC^\bullet_-(V)$.

Finally, using the acyclicity of \cbr{V} we see that the projection map
\[CC^\bullet_-(V)\longrightarrow \Sigma\widehat{CC}^\bullet_-(V)\]
is a quasi-isomorphism.  This finishes the proof.
\end{proof}

Next we shall make a link between negative cyclic cohomology of an $A_\infty$-algebra and a more customary notion of the negative cyclic \emph{homology} of a differential graded algebra. This departure from our convention to work with cohomology rather than homology is necessary because the equivariant cohomology of a loop space is expressed in \cite{jjones} in terms of cyclic homology rather than cyclic cohomology.

Recall that for a DGA $V$ its cyclic Tsygan complex is the following bicomplex $CC_\bullet^-(V)$ lying in the right half-plane and formed by taking direct products:
\[\xymatrix{ \ldots & \Cbr{V} \ar_-{1-z}[l] & \Choch{V}{V} \ar_-N[l] & \Cbr{V} \ar_-{1-z}[l] }.\]
Here \Choch{V}{V} is the usual homological Hochschild complex of $V$ and \Cbr{V} is the bar-complex of $V$ (acyclic in the unital case). The operators $1-z$ and $N$ are formed as in the cohomological cyclic complex.

Then we have the following almost obvious result.

\begin{lemma}\label{finiteness}
Let $V$ be a (unital) differential graded algebra whose Hochschild homology \Hhoch{V}{V} is finite dimensional in each degree. Then $HC^-_\bullet(V)$ is isomorphic to the graded dual of the graded $\gf$-vector space $HC^\bullet_-(V)$.
\end{lemma}

\begin{proof}
Note that the finite dimensionality assumption ensures that the complex $[\choch{V}{V^*}]^*$ is quasi-isomorphic to \Choch{V}{V}. Since the functor $?^*:=\Hom_\gf(?,\gf)$ takes direct sums to direct products we conclude that $[CC^\bullet_-(V)]^*$ is quasi-isomorphic to $CC^-_\bullet(V)$.
\end{proof}

\begin{rem}\label{simpleconnectivity}
Let $V$ be a differential graded algebra such that $H^n(V)=0$ for $n<0$, $H^0(V)=\gf$, $H^1(V)=0$ and $H^n(V)$ is finite dimensional for all $n$. For example, the cochain algebra $C^\bullet(X)$ of a simply-connected space $X$ of finite type satisfies these conditions. Then it is easy to see from the spectral sequence associated with the normalised Hochschild complex $\overline{C}^\bullet_\mathrm{Hoch}(V,V^*)$ that \hhoch{V}{V^*} is finite dimensional in each degree.
\end{rem}

With these preparations we are ready to formulate our last result.

\begin{theorem}
Let $M$ be a rational Poincar\'e duality space of formal dimension $d$. Then $H_\bullet^{S^1}(LM)$ has the structure of a graded Lie algebra of degree $2-d$. If $d\neq 0 \mod 4$ then two homotopy equivalent spaces give rise to isomorphic graded Lie algebras. If $d=0\mod 4$ then the choice of the fundamental cycle in $H_d(M)$ leads to two possibly nonisomorphic  graded Lie algebra structures on  $H_\bullet^{S^1}(LM)$.
\end{theorem}

\begin{proof}
We know from \cite{jjones} that $H^n_{S^1}(LM)\cong HC_{-n}^-(C^\bullet(M))$ and using Lemma \ref{finiteness} together with Remark \ref{simpleconnectivity} we conclude that $H_n^{S^1}(LM)\cong HC^{-n}_-(C^\bullet(M))$.

Furthermore, the differential graded algebra $C^\bullet(M)$ clearly satisfies the conditions of Lemma \ref{comparison} with $N=0$ since $\hhoch{C^\bullet(M)}{[C^\bullet(M)]^*} \cong H_{-\bullet}(LM)$ and $H_{\bullet}(LM)$ is concentrated in nonnegative degrees. Therefore we have an isomorphism
\begin{equation} \label{shift}
H_n^{S^1}(M)\cong HC^{-n}_-(C^\bullet(M))\cong HC^{-n+1}(C^\bullet(M)).
\end{equation}
Similarly to the proof of Theorem \ref{loopproduct} let $(V,m)$ be a minimal symplectic $A_\infty$-algebra weakly equivalent to $C^\bullet(M)$. Since by Proposition \ref{lem_duciso} the complex \cchoch{V} is isomorphic (with an appropriate shift) to the complex $S\choch{V}{V}$ consisting of symplectic vector fields, its homology supports the structure of a graded Lie algebra (the Gerstenhaber bracket). Since $\hchoch[n]{C^\bullet(M)} \cong \hchoch[n]{V}$ we have the Lie bracket on \hchoch[n]{C^\bullet(M)}:
\[\hchoch[n]{C^\bullet(M)} \otimes \hchoch[k]{C^\bullet(M)} \rightarrow \hchoch[n+k-d+1]{C^\bullet(M)}.\]
The latter is translated into the string bracket via isomorphism \eqref{shift}:
\[H^{S^1}_n(LM)\otimes H^{S^1}_k(LM)\rightarrow H^{S^1}_{n+k-d+2}(LM).\]

Now let us consider the problem of homotopy invariance.  For a given rational Poincar\'e duality space $N$ homotopy equivalent to $M$, the cohomology algebras $H^\bullet(M)$ and $H^\bullet(N)$ are isomorphic but the matrices of inner products on them differ by a nonzero rational number coming from the choice of the fundamental class. If $d\neq 0\mod 4$ then we see, using appropriate rescalings that $H^\bullet(M)$ and $H^\bullet(N)$ are isomorphic as Frobenius algebras after all. It follows that in this case the graded Lie algebra structures on  $H^{S^1}_\bullet(LM)$ and $H^{S^1}_\bullet(LN)$ are isomorphic. If  $d=0\mod 4$ then the signatures of $M$ and $N$ might differ by a sign depending on whether the given homotopy equivalence $M\rightarrow N$ preserves orientation or changes it.  So we have precisely two nonisomorphic structures  of a symplectic $A_\infty$-algebra on $H^\bullet(M)$ depending on the choice of orientation. This leads to two possibly nonisomorphic graded Lie algebras on $H^{S^1}_\bullet(LM)$.
\end{proof}

\begin{rem}
Clearly, the loop bracket and the loop product as defined above are nontrivial as a quick calculation of the Hochschild cohomology of $H^\bullet(S^n)$ makes clear. Since the cyclic Hochschild cohomology of a truncated polynomial algebra is concentrated in even dimensions, we conclude that for spaces such as $S^n$ and $\mathbb{C}P^n$ the string bracket is zero. This is in agreement with the calculations made in \cite{fethvi}. To obtain an example of a nontrivial string bracket consider the space $X:=S^3\times S^3$ and let $A:=H^\bullet(X)$. Then the Lie algebra consisting of symplectic vector fields of degree $0$ is identified with the Lie algebra of the group of automorphisms of $H^3(X)$ preserving the Poincar\'e duality form. In other words this is the algebra of traceless $2\times 2$ matrices and it clearly has nontrivial commutators, therefore the string bracket is nontrivial in $H_\bullet^{S^1}(X)$.
\end{rem}

\end{document}